\documentclass[a4paper,12pt]{article}
\usepackage[english]{babel}
\usepackage{amssymb}
\usepackage{amsthm}
\usepackage{amsmath}
\usepackage{graphicx}
\usepackage{array,hhline}
\usepackage[]{algorithm2e}

\newtheorem{lemma}{Lemma}[section]
\newtheorem{theorem}[lemma]{Theorem}
\newtheorem{proposition}[lemma]{Proposition}

\theoremstyle{definition}
\newtheorem{definition}[lemma]{Definition}
\newtheorem{remark}[lemma]{Remark}

\numberwithin{equation}{section} \numberwithin{figure}{section}

\newcommand{\Xset}{\mathcal{X}}

\newcommand{\upperRomannumeral}[1]{\uppercase\expandafter{\romannumeral#1}}

\begin{document}
%%%%%%%%%%%%%%%%%%%%%%%%%%%%%%%%%%%%%%%%%%%%%%%%%%%%%%%%%%%%%%%%%%%%%%%%
\title{On bivariate fundamental polynomials}

\author{V. Vardanyan {\normalsize({\tt
vahagn.vardanyan94@gmail.com})}\\
Department of Mathematics and Mechanics\\
Yerevan State University\\
A. Manukyan St. 1\\
0025 Yerevan, Armenia\\}

\date{}

\maketitle

%%%%%%%%%%%%%%%%%%%%%%%%%%%%%%%%%%%%%%%%%%%%%%%%%%%%%%%%%%%%%%%%%%%%%%%%
%%%%%%%%%%%%%%%%%%%%%%%%%%%%%%%%%%%%%%%%%%%%%%%%%%%%%%%%%%%%%%%%%%%%%%%%
\begin{abstract}
An $n$-independent set in two dimensions is a set of nodes admitting
(not necessarily unique) bivariate interpolation with polynomials of
total degree at most $n.$ For an arbitrary $n$-independent node set
$\mathcal X$ we are interested with the property that each node
possesses an $n$-fundamental polynomial in form of product of linear
or quadratic factors. In the present paper we show that each node of
$\mathcal X$ has an $n$-fundamental polynomial, which is a product
of lines, if $\#\mathcal X\le 2n+1.$ Next we prove that each node of
$\mathcal X$ has an $n$-fundamental polynomial, which is a product
of lines or conics, if $\#\mathcal X\le 2n+[n/2]+1$. We have a
counterexample in each case to show that the results are not valid
in general if $\#\mathcal X\ge 2n+2$ and $\#\mathcal X\ge
2n+[n/2]+2,$ respectively.
\end{abstract}

{\bf Key words:} Bivariate polynomial, interpolation, fundamental
polynomial, conic, $n$-poised, $n$-independent nodes.

{\bf Mathematics Subject Classification (2010):} \\
primary: 41A05, 41A63; secondary 14H50.

%%%%%%%%%%%%%%%%%%%%%%%%%%%%%%%%%%%%%%%%%%%%%%%%%%%%%%%%%%%%%%%

\vspace{5mm}

\section{Introduction}

Let $\Pi_n$ be the space of bivariate polynomials of total degree at most $n:$
\begin{equation*}
\Pi_n=\left\{\sum_{i+j\leq{n}}a_{ij}x^iy^j:a_{ij}\in \mathbb{R}
\right\}.
\end{equation*}
We have that
\begin{equation*}
N:=N_n:=\dim \Pi_n=\binom{n+2}{2}.
\end{equation*}
Consider a set of distinct nodes (points)
\begin{equation*}
\Xset_s=\{ (x_1, y_1), (x_2, y_2), \dots , (x_s, y_s) \} .
\end{equation*}
The problem of finding a polynomial $p \in \Pi_n$ which satisfies
the conditions
\begin{equation}\label{intpr}
p(x_i, y_i) = c_i, \quad i = 1, 2, \dots s  ,
\end{equation}
is called interpolation problem. A polynomial $p \in \Pi_n$ is
called an $n$-fundamental polynomial for a node $A = (x_k, y_k) \in
\Xset_s$ if
\begin{equation*}
p(x_i, y_i) = \delta _{i k}, \quad  i = 1, \dots , s ,
\end{equation*}
where $\delta$ is the Kronecker symbol. We denote this fundamental
polynomial by $p_k^\star=p_A^\star = p_{A, \Xset_s}^\star.$
Sometimes we call fundamental also a polynomial that vanishes at all
nodes of $\Xset$ but one, since it is a nonzero constant times a fundamental
polynomial.

\begin{definition}
A set of nodes $\Xset$ is called \emph{$n$-independent} if all its nodes
have fundamental polynomials. Otherwise, $\Xset$ is called
$n$-\emph{dependent.}
\end{definition}
Fundamental polynomials are linearly independent. Therefore a
necessary condition of $n$-independence is $\#\Xset \le N.$ Having
fundamental polynomials of all nodes of $\mathcal X$ we get a
solution of general interpolation problem \eqref{intpr} by using the
Lagrange formula:
\begin{equation}\label{L}
p(x, y) = \sum_{i=1}^s c_ip_i^\star(x,y).
\end{equation}
Thus we get readily that the node set $\Xset_s$ is $n$-independent
if and only if it is $n$-solvable, meaning that for any data $\{c_1,
\dots , c_s \}$ there exists a (not necessarily unique) polynomial
$p \in \Pi_n$ satisfying the conditions \eqref{intpr}.

\begin{definition}
The interpolation problem with the set of nodes $\Xset_s$ is called
$n$-\emph{poised} if for any data $\{c_1, \dots, c_s\}$ there exists a
unique polynomial $p \in \Pi_n$, satisfying the conditions
\eqref{intpr}.
\end{definition}

\noindent A necessary condition for $n$-poisedness is $s=\#\Xset_s = N.$ We have also that
a set $\Xset_N$ is $n$-poised if and only if it is $n$-independent. The following proposition
is based on an elementary Linear Algebra argument.

\begin{proposition}\label{poised}
The interpolation problem with the set of nodes $\Xset_N$ is
$n$-poised if and only if the following condition holds:
$$ p\in\Pi_n, \ p(x_i,y_i)=0, \ i=1,\ldots,N \Rightarrow p=0.$$
\end{proposition}

\noindent Now let us bring some results on $n$-independence we shall use in the sequel.
Let us start with the following simple but important result of Severi (see \cite{S}):
\begin{theorem} [\cite{S}] \label{Severi} Any set $\mathcal X,$ with
$\#{\mathcal X}\le n+1,$ is $n$-independent.
\end{theorem}
\begin{remark}  \label{rem} For each node $A\in \mathcal X$ here we can find $n$-fundamental
polynomial which is a product of $\#\mathcal X-1\le n$ lines, each of
which passes through a respective node of $\mathcal X\setminus\{A\}$
and does not pass through $A.$
\end{remark}
Next two results extend the Severi theorem to the cases of sets with no more than $2n+1$ (see \cite{E}, Proposition 1) and $3n-1$ (see \cite{HM}, Theorem 5.3) nodes, respectively.
\begin{theorem}  [\cite{E}] \label{thind2n+1} Any set $\mathcal X,$ with $\#{\mathcal
X}\le 2n+1,$ is $n$-independent, if and only if no $n+2$ nodes of
$\mathcal X$ are collinear.
\end{theorem}

\begin{theorem} [\cite{HM}] \label{th3n-1} Let $\mathcal X$ be set of nodes with $\#{\mathcal X}\le 3n.$
 Then the set $\mathcal X$ is $n$-dependent if and only if one of the following hold:\hfill\break
 i) $n+2$ nodes of $\mathcal X$ are collinear,\hfill\break
 ii) $2n+2$ nodes of $\mathcal X$ are lying on a conic, \hfill\break
 iii) $\#{\mathcal X}=3,$ there are curves $\gamma\in \Pi_3$ and $p\in \Pi_n$ such that $\gamma\cap p = {\mathcal X}.$
\end{theorem}

Here we use the same letter, say $p$, to denote the polynomial
$p\in\Pi_n\setminus \Pi_0$ and the algebraic curve defined by the
equation $p(x, y)=0.$ We denote lines and conics by $\alpha$ and
$\beta,$ respectively.

Note that, according to Theorem \ref{poised}, the
interpolation problem with $\Xset_N$ is $n$-poised if and only if
there is no algebraic curve of degree $\le n$ passing through all
the nodes of $\Xset_N.$

At the end of this section let us discuss  the problem we consider.
In view of the Lagrange formula \eqref{L} it is very important to
find $n$-independent (i.e., $n$-solvable) sets for which the
fundamental polynomials have the simplest possible forms. In Section
2 we characterize $n$-independent sets for which all fundamental
polynomials are products of lines. It is worth mentioning that for
the natural lattice, introduced by Chung and Yao in \cite{CY}, the
fundamental polynomials have the mentioned forms. But in this case
the nodes satisfy very special conditions. Namely, they are
intersection points of some $n+2$ given lines. In our
characterization (see forthcoming Theorem \ref{th2n+1}, Proposition
\ref{prop2n+1}) the restrictions on the node set are much more weak.
In Sections 3 we consider a much more involved problem. Here we
characterize $n$-independent node sets for which all fundamental
polynomials are products of lines or conics.

\section{The fundamental polynomials as products of lines}

\begin{theorem} \label{th2n+1}
\label{lm} Let $\mathcal X$ be an $n$-independent set of nodes with $\#{\mathcal
X}\le 2n+1.$
Then for each
node of $\mathcal X$ there is an $n$-fundamental polynomial, which is a
product of lines. Moreover, this statement is not true in
general for $n$-independent node sets $\mathcal X$ with $\#{\mathcal X}\ge 2n+2$
and $n\ge 2.$ \end{theorem}

The first statement of Theorem follows from the following result which covers more wider setting.
\begin{proposition} \label{prop2n+1}
\label{lm} Let $\mathcal X$ be a set of nodes with $\#{\mathcal
X}\le 2n+1$ and $A\in {\mathcal X}.$ Then the following three
statements are equivalent\hfill\break
i) The node $A$ has an $n$-fundamental polynomial,\hfill\break
ii) The node $A$ has an $n$-fundamental polynomial, which is a product of
linear factors,\hfill\break
iii) No $n+1$ nodes of $\mathcal X\setminus\{A\}$ are collinear
together with the node $A.$
 \end{proposition}

\section{The fundamental polynomials as products of lines and conics}

\begin{theorem}
\label{conm} Let $\mathcal X$ be an $n$-independent set of nodes with $\#{\mathcal
X}\le 2n+[n/2]+1.$
Then for each
node of $\mathcal X$ there is an $n$-fundamental polynomial, which is a
product of lines and conics. Moreover, this statement is not true in
general for $n$-independent node sets $\mathcal X$ with $\#{\mathcal X}\ge 2n+[n/2]+2$ and $n\ge 3.$
\end{theorem}

The first statement of Theorem follows from the following result
which covers more wider setting.
\begin{proposition}\label{himn}
\label{lm} Let $\mathcal X$ be a set of nodes with $\#{\mathcal
X}\le 2n+[n/2]+1$ and $A\in \mathcal X.$ Then the following three
statements are equivalent:\hfill\break i) The node $A$ has an
$n$-fundamental polynomial,\hfill\break ii) The node $A$ has an
$n$-fundamental polynomial, which is a product of lines and
conics,\hfill\break iii) a) no $n+1$ nodes of $\mathcal X\setminus
\{A\}$ are collinear together with $A$,\hfill\break b) if $n+1$
nodes of $\mathcal X\setminus \{A\}$ are collinear and are lying in
a line $\alpha$ then no
 $n$ nodes of $\mathcal X\setminus (A\cup\alpha)$ are collinear together with $A$,\hfill\break
 c) no $2n+1$ nodes of $\mathcal X\setminus \{A\}$ are lying on an irreducible conic together with $A$.
 \end{proposition}

\end{document}